\documentclass[12pt]{article}
\usepackage{amssymb}
\usepackage{amsthm}

\newcommand{\Z}{{\mathbb Z}}

\newtheorem{theorem}{Theorem}
\newtheorem{corollary}{Corollary}
\newtheorem{lemma}{Lemma}

\title{On equitable zero sums}

\author{Ernie Croot and Christian Elsholtz}

\begin{document}

\maketitle
\section{Introduction}
There is a rich literature on conditions guaranteeing that certain sums of 
residue classes cover certain other residue classes modulo some 
integer $N$.
Often it is of particular importance for applications to know that the 
class $0 \bmod N$ can be represented as a sum of the studied residue 
classes, and the zero class is often the most difficult case.
A well-known result along these lines is the famous 
Erd\H{o}s-Ginzburg-Ziv theorem, which says that
any sequence of $2N-1$ integers contains a subsequence of $N$ integers 
whose sum is zero; furthermore, the example of $N-1$ copies of $0$ 
and $N-1$ copies of $1$ shows that all residue classes, 
except the zero class, can be represented if the sequence were only of length
$2N-2$.

In the study of sums of distinct residue classes 
modulo $N$ (for example, the work of Olson \cite{Olson:1968}), the example
$a_1=1, \ldots , a_r=r$, with $r=  \sqrt{4N}-1$,
shows that again the zero residue class is the most difficult to 
represent.
\bigskip

Given a sequence of integers, by a {\it subsum} we mean the sum of elements of
some subsequence.   It is well-known (see Lemma \ref{Erdoslemma} below)
that if $N \geq 2$ is 
an integer, and if 
$$
a_1,...,a_M\ \in\ \Z
$$ 
is a sequence of integers with 
$$
M\ \geq\ N,
$$ 
then at least one of the subsums is $0 \pmod{N}$.  However, it can be
the case that there is just one subsequence that satisfies this, 
as is the case when
$$
M\ =\ N,\ {\rm and\ } a_1 = \cdots = a_M = 1.
$$
Nonetheless, if $M$ is large enough, say size about $2N$, then we would expect that
there are many subsequences leading to a subsum that is $0 \pmod{N}$.
Even so, it is easy to construct examples where there are much fewer than 
the expected number of such subsequences, which is $2^M/N$:  Say we take 
$M = 3N$, and again take our sequence to consist of all $1$'s.  Then, there are only
$$
{3N \choose N} + {3N \choose 2N} + {3N \choose 3N}\ =\ 2 {3N \choose N} + 1
$$
subsequences whose sum-of-elements is $0 \pmod{N}$.  

What we prove in the present paper is that if $M \geq 4N$, then 
there exists a subsequence of size at least $N$, such that this subsequence
contains at least the expected number of subsums equal to 
$0$ modulo $N$, at least when $N \geq 3$ is odd.  When this happens, we say that our
subsequence is {\it equitable}, for obvious reasons:

In the example above, chosing any 2N of the 3N terms (of $1$'s), gives a sequence
with 
$$
{2N \choose N}+1\ \gg\ 2^{2N}/\sqrt{N}
$$ 
many zero sums.  So, that 
sequence of length $2N$ (of $1$'s) is an {\it equitable} subsequence.

So, our theorem proves that there is a subsequence with a higher density 
of zero sums. 
A result of this type may be of interest in the study of arithmetic Ramsey 
theory. It may very well be the case that such results could follow from 
general reagularity results on graphs or hypergraphs, albeit the constants
involved might be very weak.

\begin{theorem}\label{main_theorem}  
For an odd integer $N \geq 3$ and any sequence of at least $4N$ integers,
there exists a subsequence of length 
$$
L\ >\ N,
$$
containing at least 
$$
2^L/N
$$
sub-subsequences whose sum-of-elements is $0$ modulo $N$.
\end{theorem}

We conjecture the following stronger version of this theorem:
\bigskip

\noindent {\bf Conjecture.}  Let $N \geq 2$ be an integer.  Then, any sequence
of integers of length at least $2N$, contains a subsequence of length $L \geq N$,
cotaining at least $2^L/N$ sub-subsequences whose sum-of-elements is $0 \pmod{N}$.
\bigskip

\section{Proof of Theorem \ref{main_theorem}}
\begin{lemma}{\label{Erdoslemma}}
If $N \geq 2$ is
an integer, and if
$$
a_1,...,a_M\ \in\ \Z
$$
is a sequence of integers with
$$
M\ \geq\ N,
$$
then at least one of the subsums is $0 \pmod{N}$.
\end{lemma}
\begin{proof}
Just take the $M\geq N$ partial sums $a_1+a_2+ \cdots + a_i$.
If there are any two sums which are congruent modulo $N$,
say $a_1+a_2+ \cdots + a_{i}$ and
$a_1+a_2+ \cdots + a_{j}$, $i<j$, then 
$a_{i+1}+a_{i+2}+ \cdots + a_{j}\equiv 0 \bmod N$.
Since $M\geq N$ it is not possible
that all partial sums are distinct modulo $N$ and not congruent to $0$.
\end{proof}

Let 
$$
a_1, ..., a_M,\ M\ \geq\ 4N
$$
be our initial sequence.  

The equitable subsequence we extract from $a_1,...,a_M$ will have the 
following properties:
\begin{itemize}
\item[a)] For each integer $n \not \equiv 0 \pmod{2N}$, the number of 
terms in the equitable subsequence 
that are congruent to $\pm n$ modulo $2N$, is even;
\item[b)]  the equitable subsequence has $L > N$ elements; and,
\item[c)]  the sum of the terms of the equitable subsequence 
is congruent to $0$ modulo 
$2N$.
\end{itemize}

Let us suppose that we can extract a subsequence of $a_1,...,a_M$ with these
properties, and after reordering terms, suppose it is 
$$
a_1,\ ...,\ a_L.
$$
Then, the number of subsequences of this new sequence that have sum $0$ modulo $N$
is given by
\begin{eqnarray*}
{1 \over N} \sum_{b=0}^{N-1} \prod_{j=1}^L \left (1 + e^{2\pi i b a_j/N} \right )
\ &=&\ {2^L \over N} \sum_{b=0}^{N-1} e^{\pi i b (a_1 + \cdots + a_L)/N}
\prod_{j=1}^L \cos(\pi b a_j/N) \\
&=&\ {2^L \over N} \sum_{b=0}^{N-1} \prod_{j=1}^L \cos(\pi b a_j/N).
\end{eqnarray*}
To obtain this last line we have used property c) above.

Now, from proprty a) above we have that for each value of $b$, this 
product of cosines is non-negative, because we may write this product as 
$$
\prod_{0 < n \leq N} \prod_{1 \leq j \leq L \atop b a_j \equiv \pm n \pmod{2N}}
\cos(\pi n/N).
$$
Note for each $n$ there are an even number of values of $j$ satisfying 
$b a_j \equiv \pm n \pmod{2N}$.

We conclude that each term in our sum over
$b$ is a non-negative real number, so the total sum is at least equal to the 
contribution of the term $b=0$, which is 
$$
2^L/N.
$$
Thus, the subsequence satisfies the conclusion of Theorem \ref{main_theorem},
and we are done, provided we can produce the sequence $a_1,...,a_L$.

\section{Using the lemma to produce the equitable subsequence}

To produce our equitable sequence, we begin by producing an auxillary sequence
$$
c_1, ..., c_T,
$$
formed by pairing up some elements $a_{i_n}, a_{j_n}$ drawn from
$a_1,...,a_M$.  This pairing is to satisfy
$$
a_{i_n}\ \equiv\ \pm a_{j_n} \pmod{2N},
$$
and is done so that $i_1,j_1,i_2,j_2,...,i_T, j_T$ are all distinct.  Note that there
may be some terms in $a_1,...,a_M$ that cannot be paired with another.

The way we define the $c_i$'s is  
$$
c_n\ :=\ a_{i_n} + a_{j_n}.
$$
This gets us the terms $c_1,...,c_V$, which are those where 
$a_{i_n}, a_{j_n} \not \equiv 0 \pmod{2N}$; the remaining terms $c_{V+1},...,c_T$ 
are all to be $0$, and correspond to sums of pairs of terms from 
$a_1,...,a_M$, both congruent to $0$ modulo $2N$. 

Noting that there can be at most 
$N$ terms from the sequence $a_1,...,a_M$ that are unpairable (because there are
at most $N$ pairs of residue classes $\pm n \not \equiv 0  \pmod{2N}$),
we deduce that
$$
2T\ \geq\ M - N\ \ \Longrightarrow\ \ T\ \geq\ 3N/2.
$$

Next, we will need to apply the following standard 
corollary of Lemma \ref{Erdoslemma} that
any sequence of length $N$ has an associated subsum that is $0 \pmod{N}$:

\begin{corollary}  Suppose that $c_1,...,c_T$ is some sequence of integers with 
$T \geq N$.  Then, there exists a subsequence of length exceeding $T-N$ 
with a subsum that is $0$ modulo $N$.
\end{corollary}

\begin{proof} The proof amounts to applying 
Lemma \ref{Erdoslemma} iteratively:  First, $c_1,...,c_N$ has a $0 
\pmod{N}$ subsum.  
Now, delete that subsequence from $c_1,...,c_T$ having sum $0 \pmod{N}$, 
and relabel terms as 
$$
c_1,...,c_U,\ {\rm where\ } T-N \leq U \leq T-1.
$$  
We then apply Lemma \ref{Erdoslemma} to {\it that} sequence, and iterate, 
until 
we reach a residual sequence having no subsum that is $0 \pmod{N}$,
and therefore has length smaller than $N$.  The union of all the terms that we deleted,
among all steps of our iterative process, has sum divisible by $N$, 
and has size at exceeding $T-N$.
\end{proof}
\bigskip

Applying this Proposition, we find that our sequence $c_1,...,c_T$ must have a 
subsequence of length exceeding
$$
T - N\ \geq\ 3N/2 - N\ =\ N/2,
$$
whose sum-of-elements is $0$ modulo $N$, which therefore means it is $0$ modulo 
$2N$, as each term $c_1,...,c_T$ is even.  But this sub-sequence of 
$c_1,...,c_T$ of length $> N/2$ having zero sum modulo $2N$ corresponds to a 
subsequence of $a_1,...,a_M$ of length $> N$.  If we write this subsequence 
(after reordering and relabeling) as 
$$
a_1,...,a_L,\ L\ >\ N,
$$
then it clearly satisfies all three properties a), b) and c) required for 
our equitable sequence, and 
so we are done.

\ \\
E.~Croot,
Georgia Institute of Technology,
School of Mathematics,
103 Skiles,
Atlanta, GA 30332, USA.\\
Email: ecroot@math.gatech.edu\\
\ \\
C.~Elsholtz,
Department of Mathematics,
Royal Holloway University of London,
Egham,
Surrey TW20 0EX,
UK.\\
Email: christian.elsholtz@rhul.ac.uk

\end{document}